\documentclass[12pt]{amsart}
\usepackage[footskip=1cm, headheight = 16pt, top=3.7cm, bottom=3cm,  right=2.5cm,  left=2.5cm ]{geometry}
\usepackage[foot]{amsaddr}
\usepackage[leqno]{amsmath}
\usepackage{lmodern}
\usepackage[T1]{fontenc}
\usepackage[utf8]{inputenc}
\usepackage[shortlabels]{enumitem}
\usepackage[center]{caption}
\usepackage{latexsym}
\usepackage{amsfonts}
\usepackage{graphicx}
\usepackage{caption}
\usepackage{xcolor}
\usepackage{todonotes}
\definecolor{favred}{RGB}{100,0,50}
\definecolor{Darkblue}{rgb}{0,0,0.4}
\definecolor{coolblack}{rgb}{0.0, 0.18, 0.39}
\usepackage{amssymb, amsthm,latexsym, mathdots, mathtools, graphicx, mathtools, xkeyval, tikz, tkz-euclide, adjustbox, float, array, subcaption, listings}
\usepackage{comment}
\usepackage[colorlinks=true,
linkcolor =blue,
citecolor=red,
urlcolor=red]{hyperref}
\usepackage{cleveref}

\newtheorem{theorem}{Theorem}[section]
\newtheorem{lemma}[theorem]{Lemma}

\newtheorem{question}{Question}
\newtheorem{observation}{Observation}

\newcounter{tbox}
\newcommand{\clm}[1]{\vspace{3mm}\refstepcounter{tbox}\noindent{\parbox{\textwidth}{(\thetbox) \emph{#1}}}\vspace{0.3mm}}
\newcommand*{\myproofname}{Proof}

\newcommand\mat{\boxminus}

\makeatletter
\renewcommand\paragraph{\@startsection{paragraph}{4}{\z@}%
  {1.25ex \@plus1ex \@minus.2ex}%
  {-1em}%
  {\normalfont\normalsize\bfseries}}
\makeatother

\newcommand{\lpt}{\mathsf{lpt}}
\newcommand{\tw}{\mathsf{tw}}

\newcommand\mx{\mathsf{max}}
\newcommand\bd{\mathsf{bd}}

\title[Hitting all longest paths in $H$-free graphs and $H$-graphs]{Hitting all longest paths in $H$-free graphs and $H$-graphs}
\date{\today}
\author[P.T. Lima]{Paloma T. de Lima$^{*,\ddagger}$}
\author[A. Nikabadi]{Amir Nikabadi$^{\ddagger}$}
\author[P. Rzążewski]{Paweł Rzążewski$^{\dagger}$}
\address{$*$ Norwegian School of Economics, Norway.}
\address{$\ddagger$ IT University of Copenhagen, Denmark.}
\address{$\dagger$ Warsaw University of Technology \& University of Warsaw, Poland.}
\address{Lima and Nikabadi acknowledge the support of the Independent Research Fund Denmark grant agreement number 2098-00012B.}
\address{Paweł Rzążewski was supported by the National Science Centre of Poland grant 2024/54/E/ST6/00094.}

\allowdisplaybreaks
\begin{document}
\maketitle

\begin{abstract}
The \textit{longest path transversal number} of a connected graph $G$, denoted by $\lpt(G)$, is the minimum size of a set of vertices of $G$ that intersects all longest paths in $G$. We present constant upper bounds for the longest path transversal number of \textit{hereditary classes of graphs}, that is, classes of graphs closed under taking induced subgraphs.

Our first main result is a structural theorem that allows us to \textit{refine} a given longest path transversal in a graph using domination properties. This has several consequences:
\begin{itemize}
\setlength{\itemsep}{1.5pt}
\item First, it implies that for every $t \in \{5,6\}$, every connected $P_t$-free graph $G$ satisfies $\lpt(G) \leq t-2$, allowing to answer a question of Gao and Shan (2019).
\item Second, it shows that every $(\textit{bull}, \textit{chair})$-free graph $G$ satisfies $\lpt(G) \leq 5$, where a \textit{bull} is the graph obtained from a four-vertex path by adding a vertex adjacent to the two middle vertices of the path, and a \textit{chair} is the graph obtained by subdividing an edge of a claw ($K_{1,3}$) exactly once.
\item Third, it implies that for every $t \in \mathbb{N}$, every connected chordal graph $G$ with no induced subgraph isomorphic to $K_t \mat \overline{K_t}$ satisfies $\lpt(G) \leq t-1$, where $K_t \mat \overline{K_t}$ is the graph obtained from a $t$-clique and an independent set of size $t$ by adding a perfect matching between them.
\end{itemize}

\noindent Our second main result provides an upper bound for the longest path transversal number in \textit{$H$-intersection graphs}. For a given graph $H$, a graph $G$ is called an \textit{$H$-graph} if there exists a subdivision $H'$ of $H$ such that $G$ is the intersection graph of a family of vertex subsets of $H'$ that each induce connected subgraphs. The concept of $H$-graphs, introduced by Bir\'o, Hujter, and Tuza, naturally captures interval graphs, circular-arc graphs, and chordal graphs, among others. Our result shows that for every connected graph $H$ with at least two vertices, there exists an integer $k = k(H)$ such that every connected $H$-graph $G$ satisfies $\lpt(G) \leq k$.
\end{abstract}

\section{Introduction}
\begin{figure}[t]
    \centering
    \includegraphics[width=0.33\textwidth]{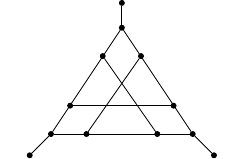}
    \caption{The Walther-Zamfirescue graph.}
    \label{fig:petersonblow}
\end{figure}

It is known that every two longest paths in a connected graph share a common vertex. An old question of Gallai~\cite{Gallai68} asks a generic form of this, whether all longest paths in a connected graph share a common vertex. According to a counterexample due to Walther~\cite{walther1969nichtexistenz} and Zamfirescu~\cite{zamfirescu1976longest}, there exists a graph (depicted in~\Cref{fig:petersonblow}) such that every vertex is omitted by some longest path of the graph, refuting Gallai's question. 
However, investigations into Gallai's question led to the formulation of a more general problem, in which one is interested in the size of the smallest set of vertices intersecting every longest path of a graph.
Let us be precise. Let $G$ be a connected graph. A \textit{longest path transversal} in $G$ is a set $S\subseteq V(G)$ such that every longest path in $G$ contains at least a vertex in $S$. The \textit{longest path transversal number} of $G$, denoted by $\lpt(G)$, is the minimum cardinality of a longest path transversal in $G$ (when we speak of a \textit{longest path}, we mean a path that achieves the maximum possible length among all paths in the graph). In this language, Gallai's question 
asked whether any graph $G$ is such that $\lpt(G)=1$. Note that the graph in \Cref{fig:petersonblow} has a longest path transversal of size two.
A graph $G$ with $\lpt(G)=3$ was constructed in~\cite{zamfirescu1976longest}. However, and surprisingly, it is not known whether there exists a connected graph $G$ with $\lpt(G)\geq4$. Accordingly, an extensive body of work has been devoted to providing non-trivial upper bounds for the longest path transversal number of a general graph~\cite{rautenbach2014transversals, long2021sublinear, kierstead2023improved}. Recently, Norin et al.~\cite{norin2025small} proved that every connected graph $G$ with $n$ vertices admits $\lpt(G) = \mathcal{O}(\sqrt{n})$, improving the previous best upper bound (see~\cite{kierstead2023improved}). 
For the special case where $G$ is chordal, Harvey and Payne~\cite{harvey2023intersecting} showed $\lpt(G) \leq 4\lceil(\omega(G))/5 \rceil$, where $\omega(G)$ is the size of the largest clique of $G$, while
Long Jr.\ et al.~\cite{long2024longest} showed that if $G$ is an $n$-vertex connected chordal graph, then $\lpt(G) = \mathcal{O}(\log^{2}n)$.

\medskip

A natural question to ask is: 

\begin{question}[see also~\cite{zamfirescu1976longest}]\label{ques:constantupper}
Is there a constant $c\in \mathbb{N}$ such that every connected graph $G$ satisfies $\lpt(G) \leq  c$?
\end{question}

Surprisingly, \Cref{ques:constantupper} remains widely open even for special classes of graphs, such as chordal graphs, though it has been resolved positively for some of its subclasses like interval graphs~\cite{balister2004longest}, dually chordal graphs~\cite{jobson2016detour}, well-partitioned chordal graphs~\cite{ahn2022} and other subclasses of chordal graphs defined by constraints on their tree intersection model~\cite{cerioli2020,long2024longest}.
Other graph classes that are known to admit constant size longest path transversals include
circular-arc graphs~\cite{balister2004longest,joos2015note} and bipartite permutation graphs~\cite{cerioli2020}.
A limited number of results also provide answers for~\Cref{ques:constantupper} when restricted to $H$-free graphs (for a given graph $H$, a graph $G$ is said to be $H$-free if it has no induced subgraph isomorphic to $H$).
 Long Jr.\ et al~\cite{long2023non} showed that if $H$ is a linear forest on at most four vertices, then every connected $H$-free graph $G$ admits $\lpt(G)=1$. 
 Gao et al.~\cite{gao2021nonempty}
 identified multiple subclasses of claw-free graphs that have $\lpt=1$. Recently, in~\cite{lima2025longest}, two of us partially improved upon the results of~\cite{long2023non} and~\cite{gao2021nonempty} by providing a complete classification of $(claw, H)$-free graphs, for $H$ of size at most five, which admit longest path transversals of size one, and by showing that $\lpt(G)=1$ when $G$ belongs to some subclasses of $P_5$-free graphs. As far as we know (and curiously enough), \Cref{ques:constantupper} is still open for the class of $P_5$-free graphs. In particular, the following question is posed in~\cite{gao2021nonempty}:

\begin{question}[see Problem 6 in~\cite{gao2021nonempty}]\label{ques:constantupper-P5}
Is there a constant $c\in \mathbb{N}$ such that every connected $P_5$-free graph $G$ admits $\lpt(G) \leq  c$?
\end{question}

\subsection*{The main results.} We prove two main results in this paper, both of which are a step toward understanding the longest path transversal numbers of graphs.
The first result, proven in~\Cref{sec:monitor-affability}, allows us to refine a given longest path transversal in a graph using domination properties. Let $G$ be a connected graph and $C$ be a component of $G$. For $M\subseteq V(G)$ let $\bd(C_i,M)$ denotes the vertices of $C_i$ with at least one neighbor in $M$. Let $t_i$ denote the length of a longest path of $G[C_i]$ with at least one endpoint in $\bd(C_i,M)$. Let $C_\mx$ be a component of $G$ such that $t_\mx\geq t_i$, $1\leq i\leq t$. We call such a component \emph{path-maximal with respect to~$M$}.

\begin{theorem}\label{thm:refining-restate}
Let $G$ be a connected graph, $M \subseteq V(G)$ be a longest path transversal of $G$ and $D\subseteq M$ be a connected dominating set of $G[M]$. Let $C_1,\dots, C_t$ be the connected components of $G\setminus M$ and let $C_\mx$ be a path-maximal component with respect to $M$. Let $S\subseteq M$ be a set that dominates $\bd(C_\mx,M)$. Then $D\cup S$ is a longest path transversal of $G$.
\end{theorem}

\Cref{thm:refining-restate} has several consequences. First, combined with the result of~\cite{chudnovsky2021}, it yields the following and thereby answers~\Cref{ques:constantupper-P5}.

\begin{theorem}\label{thm:p6-restate}
Let $t\in \{5,6\}$. Every connected $P_t$-free graph $G$ satisfies $\lpt(G) \leq t-2$.
\end{theorem}

Second, paired with the results of~\cite{hodur2025finding}, it implies the following bound for the class of $(bull,chair)$-free graphs.

\begin{theorem}\label{thm:bull-chair-restate}
Let $G$ be a connected $(bull, chair)$-free graph. Then $\lpt(G) \leq 5$.
\end{theorem}

\Cref{thm:refining-restate} also represents a step forward toward answering~\Cref{ques:constantupper} on chordal graphs. Kang et al.~\cite{kang2017width} showed that interval graphs, for which \Cref{ques:constantupper} has been answered, are ($K_3 \mat \overline{K_3}$)-free chordal graphs.
\Cref{thm:refining-restate}~implies that the longest path transversal number of chordal ($K_t\mat \overline{K_t}$)-free graphs is bounded by $t-1$.

\begin{theorem}\label{thm:chordal-lpt-lct:restate}
For every $t \in \mathbb{N}$, every connected chordal graph $G$ with no induced subgraph isomorphic to $K_t \mat \overline{K_t}$ satisfies $\lpt(G) \leq t-1$.
\end{theorem}

The second main result of this paper, proven in~\Cref{sec:h-graphs}, contributes to the longest path transversal number of $H$-intersection graphs, and hence naturally provides an upper bound for the longest path transversal number of interval graphs, circular-arc graphs, and chordal graphs, within the language of intersection graphs (see~\Cref{obs:h-graphs}). Note that, for a fixed graph $H$, our result provides a constant upper bound for the class of $H$-graphs.

\begin{theorem}\label{thm:h-graphs:restate}
For every connected graph $H$ with at least two vertices, every connected $H$-graph $G$ satisfies $\lpt(G)\leq 4(\tw(H)+1)|E(H)|$.
\end{theorem}

\newcommand{\app}{$\clubsuit$}

\section{Preliminaries}\label{sec:prelim}
\subsection*{Sets and numbers}
We use $\mathbb{N}$ to denote the set of positive integers. We let $[n]:= \{1,\dots, n\}$ for every $n\in \mathbb{N}$.
For a set $\Omega$, we denote by $2^{\Omega}$ the set of all subsets of $\Omega$. For $S \subseteq \Omega$, we let $\overline{S} = \Omega \setminus S$. Given $\mathcal{S} \subseteq 2^\Omega$, a \emph{hitting set} for $\mathcal{S}$ is a set $X \subseteq \Omega$ such that for all $S \in \mathcal{S}$, $X \cap S \neq \emptyset$.
If $X \cap S \neq \emptyset$, we may say that \emph{$X$ hits $S$}.
Let $H$ be a set. A family $\mathcal{F} = \{H_1, \dots, H_k \}$ of subsets of $H$ is said to satisfy the \textit{Helly property} if the following holds:
For every subfamily $\mathcal{F}' \subseteq \mathcal{F}$, if every pair $H_i, H_j \in \mathcal{F}'$ satisfies $H_i \cap H_j \neq \emptyset$, then $\bigcap_{ H_i \in \mathcal{F}'} \{H_i\} \neq \emptyset$.

\begin{figure}[t]
    \centering
    \includegraphics[width=0.6\textwidth]{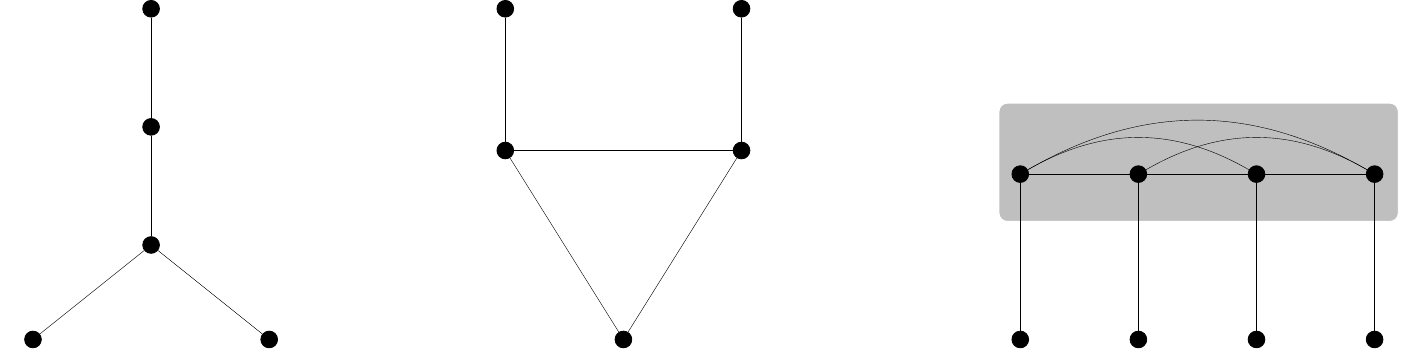}
    \caption{From left to right: a chair, a bull, and a $K_4 \mat \overline{K_4}$.}
    \label{fig:chair-bull-hedg}
\end{figure}

\subsection*{Graphs}
Graphs in this paper have finite vertex sets and no loops or parallel edges. Let $G = (V,E)$ be a graph. A \textit{clique} in $G$ is a set of pairwise adjacent
vertices. A \textit{stable set} or \textit{independent set} in $G$ is a set of pairwise non-adjacent vertices.
A subset $D \subseteq V(G)$ is a \textit{dominating set} in $G$ if each vertex of $G$ either belongs to $D$ or is adjacent to some vertex of $D$.
 For $X \subseteq V(G)$, we denote the subgraph of $G$ induced by $X$ as $G[X]$, that is $G[X] = (X, \{uv \colon u, v \in X \mbox{ and } uv \in E \})$. We denote by $N(X)$ vertices of $G\setminus X$ with a neighbor in $X$, and $N[X] = N(X)\cup X$. For disjoint $X, Y \subseteq G$, we say that $X$ is \textit{complete} to $Y$ if every vertex in $X$ is adjacent to every vertex in $Y$, and $X$ is \textit{anticomplete} to $Y$ if there are no edges between $X$ and $Y$.
We let $P = p_1p_2\dots p_k$ denote a path in $G$.
The length of $P$ is the number of edges in $P$.
We call the vertices $p_1$ and $p_k$ the endpoints of $P$ and say that $P$ is a path from $p_1$ to $p_k$. For a vertex $x\in V(G)$, we say $x$ is \textit{between $p_i$ and $p_j$}, if $x$ is in the path $p_i \dots p_j$. We also say \textit{$x$ has no neighbor from $p_i$ to $p_j$}, if there is no vertex $y$ between $p_i$ and $p_j$ such that $xy$ is an edge.

Let $G$ be a graph. A \textit{tree decomposition} of $G$ is a pair $(T, \mathcal{X})$ where $T$ is a tree and $\mathcal{X} \colon V(T) \rightarrow 2^{V(G)}$ is a
map such that:
\begin{enumerate}[\rm (i)]
    \item for every $v \in V(G)$, there exists $t \in V (T)$ such that $v \in \mathcal{X}(t)$;
    \item for every edge $uv \in E(G)$, there exists $t \in V(T)$ such that $u,v \in \mathcal{X}(t)$;
    \item for every $v \in V(G)$, the graph $T[\{t \in V(T) \colon v \in \mathcal{X}(t)\}]$ is non-empty and connected. 
\end{enumerate}
The sets $\mathcal{X}(t)$ for $t \in V(T)$ are called the \textit{bags} of $(T,\mathcal{X})$. The \textit{treewidth}
of $G$, denoted $\tw(G)$, is the minimum width of a tree decomposition of $G$, where the \textit{width} of a tree-decomposition $(T, \mathcal{X})$ is defined as $\max_{t\in V(T)} |\mathcal{X}(t)|-1$. The \textit{pathwidth} of $G$, is defined analogously with $T$ being a path instead of a tree.

For $X, Y \subseteq G$ we define the \textit{boundary of $X$ with respect to $Y$}, denoted by $\textsf{bd}(X,Y)$, as the set of vertices of $X$ that have a neighbor in $Y$. More precisely:
\[
\textsf{bd}(X,Y) = \{ x\in X \colon N_{G}(x) \cap Y \neq \emptyset\}.
\]
\subsection*{Special graphs} A graph is \emph{chordal} if it has no induced cycle of length at least four.
An \emph{interval intersection representation} for a graph $G$ is a set ${\mathcal I}=\{I_v \colon v\in V(G)\}$ of intervals, such that $uv\in E(G)$ if and only if $I_u\cap I_v\neq\emptyset$. A graph is an \emph{interval graph} if it admits an interval intersection representation.
For given graphs $G$ and $H$, we say that $G$ is \textit{$H$-free} if $G$ does not contain $H$ as an induced subgraph. We say $G$ is $(H_{1}, H_{2})$ free if $G$ does not contain $H_1$ and $H_2$ as induced subgraphs. We let $P_n$, $C_n$, and $K_n$ denote the chordless path, chordless cycle, and the complete graph on $n$ vertices. For integer $t \geq 1$, we denote by $tP_n$ the graph obtained from the disjoint union of $t$ copies of the $n$-vertex path, and for graphs $G_1$, $G_2$, we write $G_1 + G_2$ to denote the disjoint union of $G_1$ and $G_2$. A \textit{claw} is a $K_{1,3}$. A \textit{chair} is a graph obtained by subdividing an edge of a claw exactly once. A \textit{bull} is a graph obtained from a four-vertex path by adding a vertex
adjacent to the two middle vertices of the path. Let $n \in \mathbb{N}$. We denote by $K_n \mat \overline{K_n}$ the graph $\Lambda$ with
$V(\Lambda) = \{v_1^1, \dots, v_n^1, v_1^2, \dots, v_n^2 \}$ such that for all $i,j \in [n]$ the following holds (see \Cref{fig:chair-bull-hedg}):
\begin{itemize}
    \item $\{v_1^1, \dots, v_n^1\}$ is a clique, and $\{v_1^2, \dots, v_n^2\}$ is an independent set; and
    \item $v_i^1$ is adjacent $v_j^2$ if and only if $i=j$.
\end{itemize}

\begin{figure}[t]
\centering
    \includegraphics[width=0.9\textwidth]{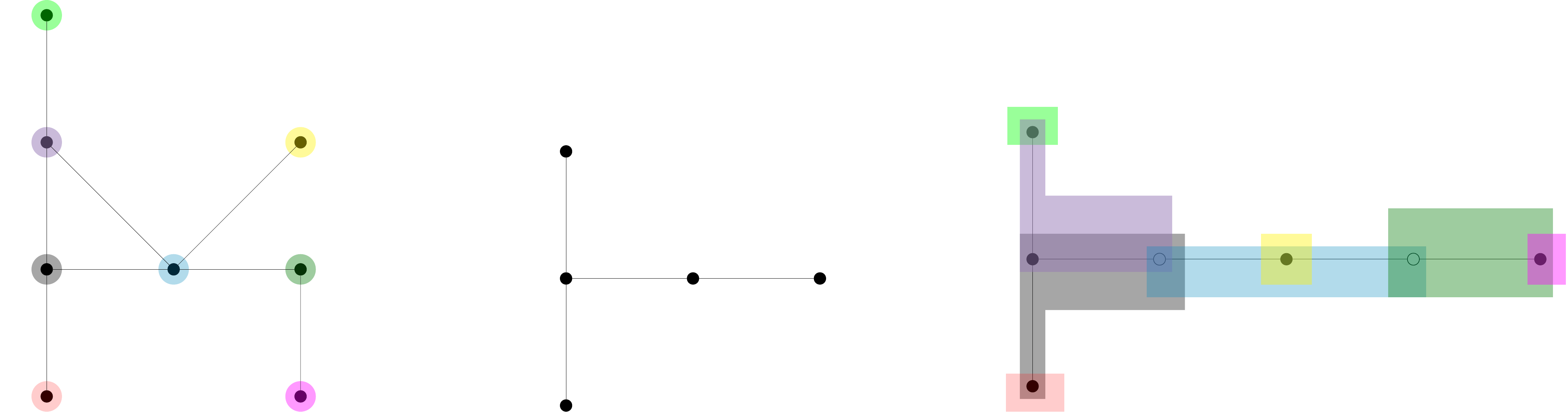}
    \caption{From left to right: A graph $G$, a graph $H$, and an $H$-representation of $G$. Vertices in white color are correspond to the subdivided edges of $H$.}
    \label{fig:h-graph-rep}
\end{figure}
\subsection*{$H$-graphs}
Let $H$ be a fixed graph. A \textit{subdivision} $H'$ of graph $H$
is a graph obtained from $H$ by replacing each edge
$e=xy$ of $H$ by a path $P_{e}$ of length at least one from $x$ to $y$ such that the interiors of the paths
are pairwise disjoint and anticomplete.
An \textit{$H$-representation} of a graph $G=(V,E)$ is a pair $(H^{\Phi}
, \Phi)$, where $H^{\Phi}$
is a subdivision of $H$ and $\Phi \colon V(G) \rightarrow 2^{V(H^{\Phi})}$ is a map such that:

\begin{enumerate}[\rm (i)]
    \item for every $v \in V(G)$, the graph $H^{\Phi}[\Phi(v)]$ is connected; and
    \item for every distinct $u, v \in V(G)$,  $uv \in E(G)$ if and only if $\Phi(u) \cap \Phi(v) \neq \emptyset$.
\end{enumerate}

\noindent We say that a graph $G$ is an \textit{$H$-graph} if $G$ admits an $H$-representation (see \Cref{fig:h-graph-rep}). We call $\Phi(v)$, $v\in V(G)$, the \textit{segments} of the representation. An $H$-graph $G$ satisfies the Helly property if the collection $\mathcal{C} = \{ \Phi(v) \colon v\in V(G) \}$ of an $H$-representation $(H^{\Phi}, \Phi)$ of $G$ satisfies the Helly property. Observe that if $H$ is a tree, then every $H$-representation satisfies the Helly property.
The concept of $H$-graphs introduced by Biró, Hujter and Tuza~\cite{biro1992precoloring}. We make the following observation, which is not difficult to check.

\begin{observation}\label{property:h-graphs}\label{obs:h-graphs}
The following hold:
\begin{enumerate}[\rm (a)]
\item Every graph $H$ is an $H$-graph.
\item $K_2$-graphs are precisely  interval graphs.
\item $K_3$-graphs are precisely  circular-arc graphs.
\item $\bigcup_{\text{tree } T} \{ T\text{-graphs}\}$ are precisely chordal graphs.
\end{enumerate}
\end{observation}
\section{Refining a transversal}\label{sec:monitor-affability}

In this section, we first prove \Cref{thm:refining-restate} which refines a given longest path transversal using domination properties. We will then apply \Cref{thm:refining-restate} and conclude that: (i) $P_6$-free graphs admit $\lpt \leq 4$, which answers the open question posed by~\cite{gao2021nonempty}; (ii) $(bull,chair)$-free graphs admit $\lpt \leq 5$, and (iii) chordal ($K_t \mat \overline{K_t}$)-free graphs admit $\lpt \leq t-1$ which adds to the extensive list of subclasses of chordal graphs that have been shown to admit constant-size longest path transversals.

We first begin with the following terminology. Given a subset $M \subset V(G)$ and the connected components $C_1,\dots C_t$ of $G\setminus M$, we let $t_i$ be the length of a longest path of $G[C_i]$ with at least one endpoint in $\bd(C_i,M)$. Recall that $\bd(C_i,M)$ denotes the vertices of $C_i$ with at least one neighbor in $M$. We let $C_\mx$ be a component such that $t_\mx\geq t_i$ for $1\leq i\leq t$. We call such a component \emph{path-maximal with respect to $M$}.

\medskip

We can now restate and prove~\Cref{thm:refining-restate}.

\begin{theorem}\label{thm:refining}
Let $G$ be a connected graph, $M \subseteq V(G)$ be a longest path transversal of $G$ and $D\subseteq M$ be a connected dominating set of $G[M]$. Let $C_1,\dots, C_t$ be the connected components of $G\setminus M$ and let $C_\mx$ be a path-maximal component with respect to $M$. Let $S\subseteq M$ be a set that dominates $\bd(C_\mx,M)$. Then $D\cup S$ is a longest path transversal of $G$.
\end{theorem}

\begin{proof}
    Let $P$ be a longest path of $G$. We want to show $V(P)\cap (D\cup S)\neq\emptyset$. 
    
    Suppose for a contradiction $V(P)\cap (D\cup S)=\emptyset$. First note that, since $M$ is a longest path transversal of $G$, $V(P) \cap M\neq\emptyset$. In particular, this implies $V(P)$ cannot be contained in $C_i$, for $1\leq i\leq t$. We will show that $V(P)\cap C_\mx=\emptyset$. We already argued $V(P)$ cannot be contained in $C_\mx$. Suppose there exists an edge $cc'\in E(P)$ such that $c\in C_\mx$ and $c'\in \bd(M,C_\mx)$. Recall that $S\subseteq M$ is a set that dominates $\bd(C_\mx,M)$. Let $s\in S$ be a vertex that dominates $c$. Since $s\in M$, let $d\in D$ be a vertex that dominates $s$. Note that it can be that $s=d$, as $S$ and $D$ do not need to be disjoint. Let $d'\in D$ be a vertex that dominates $c'$. It can also be $c'=d'$. We are now able to replace the edge $cc'$ in $P$ by the path $csd\dots d'c'$, where $d\dots d'$ is a path in $G[D]$ connecting $d$ and $d'$. Such a path is certain to exist as $D$ is connected. With this replacement, we obtained a path that is longer than $P$, as at least one vertex was added to it, a contradiction. Hence, $V(P)\cap C_\mx=\emptyset$.

    Note that $P$ cannot have an endpoint in $M$. Indeed, since $D$ dominates $G[M]$, we would be able to extend $P$ to $D$, a contradiction. Let $p_z$ be a vertex of $P$ such that $p_z\in M$, but $p_{z+j}\notin M$, for $j\geq 1$. Such a vertex exists since $p_k\notin M$. Moreover, since the $C_i$'s are connected components of $G\setminus M$, the vertices $p_{z+j}$, for $j\geq 1$, all belong to the same component. Let $C_i$ be this component. Since $p_{z+1}\dots p_k$ is a path in $G[C_i]$, its length is at most $t_i$.
    To reach a contradiction, we will modify the path $P$ in order to obtain an even longer path.
    Let $q_1\dots q_{t_1}$ be a longest path of $G[C_\mx]$ with at least one endpoint in $\bd(C_\mx, M)$. Assume $q_1\in \bd(C_\mx, M)$.
    Let $s\in S$ be a vertex that dominates $q_1$, and let $d\in D$ be a vertex that dominates $s$. Let $d''\in D$ be a vertex that dominates $p_z$.
    Consider the path $P'=p_1p_2\dots p_zd''\dots dsq_1\dots q_{t_1}$, where $d''\dots d$ is a path in $G[D]$ connecting $d''$ and $x$, which exists since $D$ is connected.
    Since $t_\mx\geq t_i$ and the length of the path $d''\dots ds$ is at least one, we conclude that $P'$ is longer than $P$, a contradiction. We can conclude $V(P)\cap (D\cup S)\neq\emptyset$, and hence $D\cup S$ is a longest path transversal of $G$. This completes the proof of \Cref{thm:refining}
\end{proof}

Before we apply~\Cref{thm:refining}, we state the following Observation of~\cite{cerioli2019}, that will be useful in the remainder of this section.

\begin{observation}[Cerioli and Lima~\cite{cerioli2019}]\label{obs:cds-lpt}
    Let $D$ be a connected dominating set of $G$. Then $D$ is a longest path transversal of $G$.
\end{observation}

We first consider $P_5$-free graphs and $P_6$-free graphs. Following~\cite{chudnovsky2021}, we say that $M \subset V(G)$ is a $monitor$ in $G$ if for every connected component $C$ of $G\setminus M$ there exists a vertex $w \in M$ that is complete to $C$. Note that a monitor is a dominating set of $G$. 

The following result of \cite{chudnovsky2021} shows that every connected $P_{t}$-free graph, $t\in \{4,5,6\}$, admits a monitor formed by the closed neighborhood of a path. In particular, this implies that this monitor itself induces a connected subgraph of $G$ and admits a small connected dominating set. More precisely:

\begin{lemma}[Chudnovsky, King, Pilipczuk, Rz{\k{a}}{\.z}ewski, and Spirkl; see Lemma 5 in~\cite{chudnovsky2021}]\label{lem:monitor}
Let $t \in \{4, 5, 6\}$, $G$ be a connected $P_6$-free graph, and $u \in V(G)$ be a vertex such that in $G$ there
is no induced $P_t$ with $u$ being one of the endpoints. Then there exists a subset $X$ of vertices such that $u \in X$,
$|X| \leq t-3$, $G[X]$ is a path whose one endpoint is $u$, and $N_{G}[X]$ is a monitor in $G$.
\end{lemma}

\begin{figure}[t]
    \centering
    \includegraphics[width=0.4\textwidth]{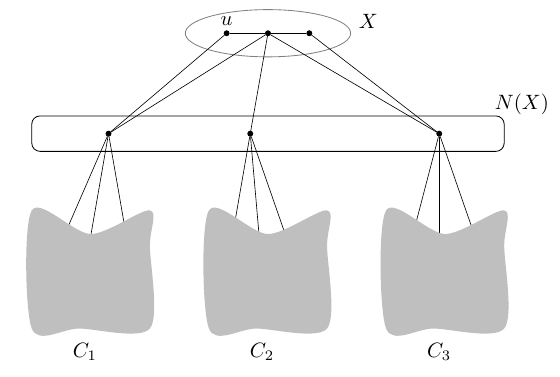}
    \caption{Outcome of \Cref{lem:monitor} for $t=6$. For $i\in [3]$, $C_i$'s are connected components in $G\setminus N_{G}[X]$.}
    \label{fig:monitor}
\end{figure}

We are now able to prove the following.

\begin{theorem}\label{thm:p6}
Let $t\in \{5,6\}$. Every connected $P_t$-free graph $G$ satisfies $\lpt(G) \leq t-2$.
\end{theorem}

\begin{proof}
    Since $G$ is $P_t$-free, for $t \in \{5,6\}$, it follows from \Cref{lem:monitor} that $G$ admits
    a monitor formed by the closed neighborhood of an induced path $X$ of length at most $t-3$. Since $N[X]$ is a monitor and $X$ is connected, $N[X]$ is also a connected dominating set of $G$, and by~\Cref{obs:cds-lpt}, it is a longest path transversal of $G$. We can now apply~\Cref{thm:refining} to refine this transversal. Let $w\in N[X]$ be the vertex that dominates a path-maximal component $C_\mx$ of $G$ with respect to $N[X]$. Since $X$ is a connected dominating set of $N[X]$ and $w$ is a vertex that dominates $\bd(C_\mx,N[X])$, by~\Cref{thm:refining}, $X\cup\{w\}$ is a longest path transversal of $G$ of size at most $t-2$. This completes the proof of~\Cref{thm:p6}.
\end{proof}

We now consider the class of $(bull,chair)$-free graphs, where monitors also play a relevant role. Indeed, the lemma below states that, in this graph class, the closed neighborhood of a maximal induced path is a monitor.

 \begin{lemma}[Hodur, Pilśniak, Prorok, and Rzążewski; see 3.2 in~\cite{hodur2025finding}]\label{lem:maximalP6inbull-chair-free}
    Let $G$ be a connected (bull, chair)-free graph on at least three vertices. Let $X\subseteq V(G)$ be such that $G[X]$ is a maximal induced path on at least three vertices. Then $N[X]$ is a monitor in $G$.
\end{lemma}

Similarly as in the proof of~\Cref{thm:p6}, we want to apply~\Cref{lem:maximalP6inbull-chair-free} to the monitor $N[X]$. The problem, however, is that~\Cref{lem:maximalP6inbull-chair-free} does not provide any bound on the length of the maximal induced path, which is crucial to obtain a constant upper bound on $\lpt$. To deal with that, we will use the following characterization of $(bull,chair)$-free graphs that have a maximal induced path with at least seven vertices.

Before we proceed, we need the following definition. A \textit{fat path} (or a \textit{fat cycle}) is a graph whose vertex set can be partitioned into nonempty sets $V_1, V_2, \dots, V_r$, such
that the sets $V_i, V_j$ are complete to each other if $j = i+1$ (or $j = i + 1 \mod r$), and anticomplete otherwise.

\begin{lemma}[Hodur, Pilśniak, Prorok, and Rzążewski; see 3.1 in~\cite{hodur2025finding}]\label{lem:RDT}
Let $G$ be a connected $(bull,chair)$-free graph, containing a maximal induced path $P$ with at least 7
vertices. If there exists a vertex $x \in V(G)$ such that $x$ is adjacent to (only) both endpoints of $P$, then we let $Q$ to be the induced cycle formed by vertices of $P$ and $x$. Otherwise, we let $Q = P$. In polynomial time one can compute a partition of $V(G)$ into sets $R, D, T$ with the following properties:
\begin{enumerate}[\rm (i)]
    \item $G[R]$ is a fat path on at least 7 vertices (or a fat cycle on at least 8 vertices);\label{propt:rdt-i}
    \item $D$ is complete to $R$, and separates $R$ and $T$;\label{propt:rdt-ii}
    \item Every component of $G - (R \cup D)$ is contained in one component of $G - N_{G}[Q]$\label{propt:rdt-iii}; and
    \item $N[R]$ is a monitor in $G$\label{propt:rdt-iv}.
\end{enumerate}
\end{lemma}

We can now consider three separate cases: when $G$ has no induced $P_6$ (in which case we simply apply~\Cref{thm:p6}), when the graph has a maximal induced path on at least seven vertices (\Cref{lem:bull-chair-maximalP7}), and when the graph has only has a maximal induced path on six vertices (\Cref{lem:bull-char-maximalP6}).

\begin{lemma}\label{lem:bull-chair-maximalP7}
 Let $G$ be a connected $(bull, chair)$-free graph. If $G$ has a maximal induced path on 
 $\ell \geq 7$ vertices, then $\lpt(G)\leq 2$.   
\end{lemma}

\begin{proof}
    By~\Cref{lem:RDT}, there is a partition of $V(G)$ into sets $R, D, T$ satisfying properties \ref{propt:rdt-i}-\ref{propt:rdt-iv}. By property \ref{propt:rdt-iv}, $N[R]$ is a monitor, and hence it is a dominating set of $G$.
    Since $R$ and $D$ are complete to each other, $N[R]$ is connected. By ~\Cref{obs:cds-lpt}, $N[R]$ is a longest path transversal of $G$.  
    Let $w\in D$ be a vertex that dominates a path-maximal component $C_\mx$ of $G$ with respect to $N[R]$. Note that $w$ dominates all vertices of $R$. Let $x$ be a vertex of $R$. Note that $\{w,x\}$ is a connected dominating set of $G[N[R]]$. Moreover, $w$ dominates $C_\mx = \bd(C_\mx,N[R])$. By~\Cref{thm:refining}, $\{w,x\}$ is a longest path transversal of $G$. This completes the proof of~\Cref{lem:bull-chair-maximalP7}.
\end{proof}

\begin{lemma}\label{lem:bull-char-maximalP6}
 Let $G$ be a connected $(bull, chair)$-free graph. If any maximal induced path of $G$ has six vertices, then $\lpt(G)\leq 5$. 
\end{lemma}

\begin{proof}
    Let $Q = q_1 \dots q_6$, be a maximal induced path of $G$. By~\Cref{lem:maximalP6inbull-chair-free}, $N[Q]$ is a monitor in~$G$, and hence also a connected dominating set of $G$. By~\Cref{obs:cds-lpt}, $N[Q]$ is a longest path transversal of $G$. Since $N[Q]$ is a monitor, let $w\in N[Q]$ be the vertex that dominates a path-maximal component $C_\mx$ of $G$ with respect to $N[Q]$. 
    
    In order to refine $N[Q]$ using~\Cref{thm:refining}, we will find a small connected set dominating $N[Q]$. 
    Let $\widetilde{Q}=Q \setminus \{q_1,q_6\}$.
    More precisely, we now show $D=\widetilde{Q}\cup \{w\}$ is a set that is connected and that dominates $N[Q]$.  We first show $w$ is complete to $Q$.
    Suppose not. If $w$ is adjacent to two consecutive vertices of $Q$, but not complete to $Q$, then there is always a set of three consecutive vertices $q_i, q_{i+1}, q_{i+2} \in Q$ such that $wq_{i+1}$ is an edge and exactly one of $q_i, q_{i+2}$ is a neighbor of~$w$. Then, the set $\{c, w, q_i, q_{i+1}, q_{i+2} \}$ induces a bull in $G$ for some vertex $c \in C_\mx$ (recall that $w$ dominates $C_\mx$).
    If $w$ has no consecutive neighbors in $Q$ but is adjacent to a vertex $q_i \in Q$ with $2\leq i\leq 5$, then the set $\{w,c,q_{i-1},q_i,q_{i+1}\}$ induces a chair in $G$.
    Since $Q$ is maximal, $w$ cannot be adjacent only to $q_1$ nor only to $q_6$. If $w$ is adjacent to both $q_1$ and $q_6$ (and only these two), then $\{q_5, q_6, w, c, q_1 \}$ induces a chair in $G$, a contradiction.
    Hence, $w$ is complete to $Q$, and $G[D]$ is connected.
    Let $v$ be a vertex of $N[Q]\setminus D$. If $v=q_1$ or $v=q_6$, then $v$ is clearly dominated by~$D$. If $v$ has a neighbor in $\widetilde{Q}$ or is adjacent to $w$, then it is also dominated by $D$.
    Otherwise, note that since $Q$ is maximal, $v$ cannot be dominated only by $q_1$ nor only by $q_6$. So $v$ must be adjacent to both $q_1$ and $q_6$ (and only these). If $v$ has a neighbor $u\notin N[Q]$, then $q_2q_3q_4q_5q_6vu$ is an induced path of length at least seven, a contradiction to the assumption that every maximal induced path of $G$ has length six. Let $c\in C_\mx$. Note that $cv\notin E(G)$, as $v$ has no neighbor outside $N[Q]$. Then $\{q_1,q_2,w,c,v\}$ induces a bull in $G$, a contradiction. Hence $D$ is indeed a connected dominating set of $N[Q]$ and, by~\Cref{thm:refining}, $D$ is a longest path transversal of $G$. This completes the proof of~\Cref{lem:bull-char-maximalP6}.
\end{proof}

\begin{figure}[t]
    \centering
  \includegraphics[width=0.55\textwidth]{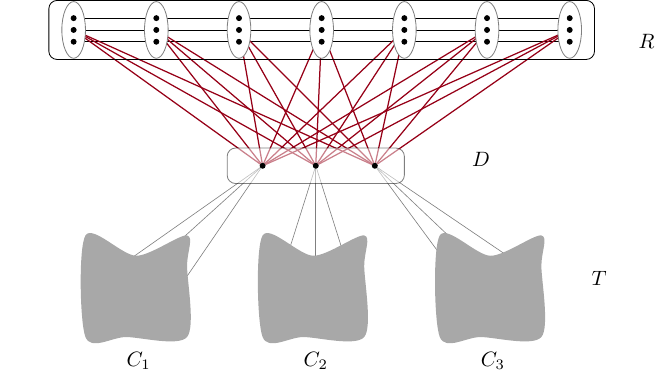}
    \caption{An $(R,D,T)$ decomposition. $C_i$'s are components in $G - (R\cup D)$. Red lines represent complete adjacency between $R$ and $D$.}
    \label{fig:(RDT)}
\end{figure}

We are now able to conclude the following.

\begin{theorem}\label{thm:bull-chair}
 Let $G$ be a connected $(bull, chair)$-free graph. Then $\lpt(G)\leq 5$.   
\end{theorem}

\begin{proof}
If $G$ is $P_6$-free, then from~\Cref{thm:p6-restate} we have $\lpt(G)\leq 4$ and the claim follows. Therefore, we may assume that there exists a set $X \subseteq V(G)$ such that $G[X]$ is an induced path on six vertices. Consider a maximal induced path $Q$ of $G$ that contains $G[X]$. If $Q$ has length at least 7, by~\Cref{lem:bull-chair-maximalP7}, $\lpt(G)\leq 2$. Otherwise, by~\Cref{lem:bull-char-maximalP6}, $\lpt(G)\leq 5$. This completes the proof of~\Cref{thm:bull-chair}.
\end{proof}

As a last application of~\Cref{thm:refining}, we now consider chordal ($K_t \mat \overline{K_t}$)-free graphs. We will make use of the observation below.

\begin{observation}[Balister, Gy\"{o}ri, Lehel, and Schelp~\cite{balister2004longest}]\label{lem:dominating-clique-inchordal}
If $G$ is a connected chordal graph, then there exists a clique $K$ in $G$ such that $K$ is a longest path transversal of $G$.
\end{observation}


\begin{theorem}\label{thm:chordal-lpt}
For every $t\in \mathbb{N}$, if $G$ is a connected chordal graph with no induced subgraph isomorphic to $K_t \mat \overline{K_t}$, then $\lpt(G) \leq t-1$.
\end{theorem}

\begin{proof}
    Let $G$ be a chordal ($K_t \mat \overline{K_t}$)-free graph. By~\Cref{lem:dominating-clique-inchordal}, there exists a clique $K$ that is a longest path transversal of $G$. Let $C_\mx$ be a path-maximal component of $G\setminus K$. Let $D\subseteq K$ be a minimal set with the property that it dominates $\bd(C_\mx,K)$. Note that $D$ also dominates~$K$, as $K$ is a clique. By~\Cref{thm:refining}, $D$ is a longest path transversal of $G$. We now show that since $G$ is ($K_t \mat \overline{K_t}$)-free, then $|D|\leq t-1$.
    
    Indeed, since $D$ is minimal, for every $d\in D$ there exists $d' \in \textsf{bd}(C_\mx,K)$ such that $d'$ is a \textit{private neighbor} of $d$, that is, $d'$ is dominated \emph{only} by $d$ in $D$. Let $\{d'_1, \dots d'_r \}$ be the set of private neighbors of the vertices $d_1, \dots, d_r$ in $D$.
    Since $G[D]$ is a clique, there are no edges between vertices in $\{d'_1, \dots d'_r \}$, otherwise, there would be an induced cycle of length four in $G$, which violates the assumption that $G$ is chordal. Thus, $ \{d_1, \dots, d_r\} \cup \{d'_1, \dots, d'_r\}$
    induces a $K_r \mat \overline{K_r}$ in $G$. Since $G$ has no induced subgraph isomorphic to $K_t \mat \overline{K_t}$, we conclude that $r \leq t-1$. This completes the proof of~\Cref{thm:chordal-lpt}.
\end{proof}

\section{$H$-graphs}\label{sec:h-graphs}
Before proving \Cref{thm:h-graphs:restate} in this section, we
introduce some notation that will be used later.

\paragraph{$H$-representation of graphs.}
Let $(H^{\Phi}, \Phi)$ be an $H$-representation of a graph $G$ (see~\Cref{sec:prelim}).
We say that $(H^{\Phi}, \Phi)$ is a \textit{nice} $H$-representation of $G$ if the following hold:
\begin{enumerate}[\rm (a)]
   \item for every $x \in  V(H^{\Phi})$,
there is a vertex $v \in  V(G)$ such that $x \in  \Phi(v)$; and
\item for every edge $xy \in  E(H^{\Phi})$,
there is $v \in V(G)$ such that $x,y \in  \Phi(v)$.
\end{enumerate}
\noindent It is known that an $H$-representation of $G$ can be turned into a nice one~\cite{chaplick2021kernelization}.
Henceforth, whenever we talk about an $H$-representation of $G$, we may assume that it is a nice one.

For a vertex $x$ of $H^{\Phi}$, let $V_x$ denote the set $\{v \in V(G) \colon x \in \Phi(v)\}$.
Note that for each $x \in V(H^{\Phi})$, the set $V_x$ is a clique in $G$.

Consider an edge $ab \in E(H)$ and fix its orientation from $a$ to $b$.
Note that the edge $ab$ of $H$ corresponds to a path in $H^{\Phi}$, denote the consecutive vertices of this path by $s_1,\ldots,s_\ell$, where $a=s_1$ and $b=s_{\ell}$.
For a vertex $v \in V_a$, let $\textsf{reach}_{a \to b}(v)$ be the maximum index $i$ such that $s_1,\ldots,s_i \in \Phi(v)$.
We emphasize that the intersection of $\Phi(v)$ with $\{s_1,\ldots,s_\ell\}$ does not have to be connected, so $\Phi(v)$ may contain some $s_j$ for $j > \textsf{reach}_{a \to b}(v)$.

\paragraph{Tree decompositions of $H$.}
For an $H$-graph $G$, by an \textit{$H$-profile of $G$} we denote a tuple $((H^{\Phi}, \Phi); (T,\mathcal{X}))$,
where $(H^{\Phi}, \Phi)$ is a nice $H$-representation of $G$ and $(T, \mathcal{X})$ is a tree decomposition of $H^{\Phi}$.
Note that since $H^{\Phi}$ is a subdivision of $H$, we have $\tw(H^{\Phi}) = \tw(H)$.

By the properties of a tree decomposition, for each $u \in V(H^{\Phi})$, the set $\{ t \in V(T) \colon u \in \mathcal{X}(t)\}$ induces a subtree of $T$. As, for every $v\in V(G)$, the set $\Phi(v)$ induces a connected subgraph of $H^{\Phi}$, we conclude that the set
$\{t \in V(T) \colon \Phi(v) \cap \mathcal{X}(t) \neq \emptyset \}$ corresponds to a subtree of $T$.
Denote this subtree by $T_{\Phi(v)}$.
For a path $P$ in $G$, we write $T_{\Phi(P)}$ as a shorthand for $\bigcup_{v \in V(P)}T_{\Phi(v)}$.

\medskip

The following lemma guarantees that, given an $H$-profile $((H^{\Phi}, \Phi); (T,\mathcal{X}))$ of a connected graph $G$, one can always find a \textit{bag} in $(T,\mathcal{X})$ such that those vertices of $G$ represented by the vertices of the bag form a longest path transversal set in $G$. More precisely:

\begin{lemma}\label{lem:h-graphs:hellytype}
Let $G$ be a connected graph and $((H^{\Phi}, \Phi); (T,\mathcal{X}))$ be an $H$-profile of $G$.
There exists $t\in V(T)$ such that
the set $\bigcup_{x \in \mathcal{X}(t)} V_x$ is a longest path transversal of $G$.
\end{lemma}

\begin{proof}
Let $\mathcal{F}$ be a family of all longest paths in $G$.
Since $G$ is connected, for every pair $P_1, P_2 \in \mathcal{F}$, there exists a vertex $x \in V(P_1) \cap V(P_2)$.
Since the set $\{T_{\Phi(P)} \colon P\in \mathcal{F} \}$ has the Helly property (as the subtrees of a tree admit the Helly property),
it follows that there exists a vertex $t \in V(T)$ such that for every $P \in \mathcal{F}$, we have $t \in T_{\Phi(P)}$.
This implies that the set $\{ v\in V(G) \colon t\in V(T_{\Phi(v)}) \} = \bigcup_{x \in \mathcal{X}(t)} V_x$ forms a longest path transversal in $G$.
\end{proof}

\paragraph{Intermediate graph.}
Let $H$ be a connected graph with at least two edges.
Take an $H$-profile $((H^{\Phi}, \Phi); (T,\mathcal{X}))$ of a connected graph $G$, where $(T,\mathcal{X})$ is a tree decomposition of $H^{\Phi}$ of width $\tw(H^{\Phi}) = \tw(H)$.
Let $t$ be given by~\Cref{lem:h-graphs:hellytype}.
We aim to show that from $\bigcup_{x \in \mathcal{X}(t)} V_x$ we might pick a subset of size bounded in terms of $H$ only,
that is still a longest path transversal of $G$.
Here we stumble on the first technical issue.
Note that the set $\mathcal{X}(t)$ might contain some vertices that are not in $H$.
On the other hand, the size of $H^{\Phi}$ can be arbitrarily large, not necessarily bounded in terms of $H$.
For this reason, we will work neither with $H$ nor with $H^{\Phi}$, but rather the \textit{intermediate graph} $\widehat{H}$.
We obtain it from $H^{\Phi}$ by contracting all edges with at least one endpoint that is not in $V(H) \cup \mathcal{X}(t)$.
Consequently, we have $V(\widehat{H}) = V(H) \cup \mathcal{X}(t)$ and thus $|V(\widehat{H})| \leq |V(H)| + \tw(H)+1 \leq 2|V(H)|$.
On the other hand, every vertex of $\mathcal{X}$ adds at most one new edge to $E(\widehat{H}) \setminus E(H)$, and thus 
\begin{equation}
\tag{S1}
    |E(\widehat{H})| \leq |E(H)|+ \tw(H)+1 \leq 2|E(H)| \label{eq:edges-of-intermediate-graph}
\end{equation}
(here we use that $|E(H)| \geq 2$).
Observe that $H^{\Phi}$ is a subdivision of $\widehat{H}$ and thus $G$ is an $\widehat{H}$-graph with  $\widehat{H}$-profile  $((H^{\Phi}, \Phi); (T,\mathcal{X}))$.
Note that for each $v \in \bigcup_{x \in \mathcal{X}(t)} V_x$, where $t\in V(T)$ is an outcome of applying~\Cref{lem:h-graphs:hellytype} on $G$, there is a vertex $x$ of $\widehat{H}$ such that $x \in \Phi(v)$. This property will be useful later.

\medskip

We can now restate and prove \Cref{thm:h-graphs:restate}.

\begin{theorem}\label{thm:h-graphs}
For every connected graph $H$ with at least two vertices, every connected $H$-graph $G$ satisfies $\lpt(G)\leq 4(\tw(H)+1)|E(H)|$.
\end{theorem}

\begin{proof}
Let $((H^{\Phi}, \Phi); (T,\mathcal{X}))$ be an $H$-profile of $G$.
Let $t\in V(T)$ be the outcome of applying~\Cref{lem:h-graphs:hellytype} on $G$, and let $\widehat{H}$ be the intermediate graph.
For simplicity of notation, we let $X=\mathcal{X}(t)$ and denote $S = \bigcup_{x \in X} V_x$.
Note that $|X| \leq \tw(H) + 1$ and each set $V_x$, for $x \in X$, is a clique in $G$, thus $S$ is partitioned into $X$ cliques.
Furthermore, recall that $S$ is a longest path transversal in $G$.

For $x \in X$, let $A_x \subseteq V(\widehat{H})$ be the set of all vertices $a$ of $\widehat{H}$ for which $V_x \cap V_a \neq \emptyset$.
For any edge $ab \in E(\widehat{H})$, where $a \in A_x$,
let $s_{x,a \to b}$ denote a vertex $v \in V_x \cap V_a$ with maximum value of $\textsf{reach}_{a \to b}(v)$; the ties are resolved arbitrarily.
Let $Q$ be the set consisting of all vertices selected this way, that is,
\[
Q = \bigcup_{x \in X} \bigcup_{a \in A_x} \bigcup_{ab \in E(\widehat{H})} \{ s_{x,a \to b} \}.
\]
We observe that since $\widehat{H}$ is connected and has at least two vertices, for each $x \in X$, the set $A_x$ is non-empty.
Consequently, for each $x \in X$, the set $Q$ intersects $V_x$.

Clearly we have $Q \subseteq S$ and, using \eqref{eq:edges-of-intermediate-graph} we get
\begin{equation}
\tag{S2}
|Q| \leq |X| \cdot 2 \cdot |E(\widehat{H})| \leq 4(\tw(H) +1)|E(H)|.\label{eq:size-of-lp}
\end{equation}
We claim the following.

\clm{For every $x \in X$, the set $Q \cap V_x$ dominates $\textsf{bd}(V(G) \setminus S,S \cap V_x)$.}\label{clm:h-graph:domset}

\medskip

For contradiction, suppose there is $x \in V$ and a vertex $u \notin S$ adjacent to some $s \in S \cap V_x$ but non-adjacent to every vertex of $Q \cap V_x$. Since $su \in E(G)$, there is some $z \in V(H^{\Phi})$ such that $z \in \Phi(u) \cap \Phi(s)$.
Let $ab \in E(\widehat{H})$ be an edge such that $z$ belongs to the path $P_{ab}$ in $H^{\Phi}$ obtained by subdividing the edge $ab$.
Note that $ab$ might not be unique if $z$ is a vertex of $V(\widehat{H})$.
Since $\Phi(s)$ is a connected subgraph of $H^{\Phi}$ and $\Phi(s)$ contains a vertex of $V(\widehat{H})$,
we conclude that $\Phi(s)$ contains the whole subpath of $P_{ab}$ from $a$ to $z$ or the whole subpath of $P_{ab}$ from $z$ to $b$.
By symmetry, assume the former. As $s \in V_x \cap V_a$, the vertex $s_{x, a\to b}$ is defined and belongs to $Q \cap V_x$.
Since $\Phi(s)$ contains the whole subpath of $P_{ab}$ from $a$ to $z$ and $s_{x, a\to b}$ is a vertex from $V_x \cap V_a$ with the maximum value of $\textsf{reach}_{a \to b}$, we observe that $z \in \Phi(s_{x, a\to b})$.
Consequently, $us_{x, a\to b} \in E(G)$, a contradiction. This proves (\ref{clm:h-graph:domset}).

\clm{$Q$ is a longest path transversal of $G$.}\label{clm:h-graph:hittingset}

\medskip

For contradiction, suppose there is a longest path $P = p_1\dots p_k$ in $G$ that avoids $Q$, that is, $V(P) \cap Q = \emptyset$.
By \Cref{lem:h-graphs:hellytype}, $S$ is a longest path transversal in $G$.
Let $i$ be the minimum index such that $p_i \in S$ and let $x \in X$ be such that $p_i \in V_x$. First, suppose that $i=1$.
Let $q_x$ be any vertex from $V_x \cap Q$; it exists as $V_x \cap Q$ is non-empty.
As $V_x$ is a clique, we obtain that $q_x,p_1\dots p_k$ is a path in $G$, contradicting the maximality of $P$. So suppose that $i > 1$.
Since $p_{i-1} \in \textsf{bd}(V(G) \setminus S,S \cap V_x)$, by \eqref{clm:h-graph:domset} there is some $q_x \in Q \cap V_x$,
such that $p_{i-1}q_x \in E(G)$.
Again using that $V_x$ is a clique, we obtain that $p_1,\ldots,p_{i-1},q_x,p_i\dots p_k$ is a path in $G$,
contradicting the maximality of $P$. This proves (\ref{clm:h-graph:hittingset}).

\medskip

From (\ref{clm:h-graph:hittingset}) and (\ref{eq:size-of-lp}) we conclude the proof of~\Cref{thm:h-graphs}.
\end{proof}


\bibliographystyle{abbrv}
\bibliography{ref}

\begin{thebibliography}{10}

\bibitem{ahn2022}
J.~Ahn, L.~Jaffke, O.~joung Kwon, and P.~T. Lima.
\newblock Well-partitioned chordal graphs.
\newblock {\em Discrete Mathematics}, 345(10):112985, 2022.

\bibitem{balister2004longest}
P.~N. Balister, E.~Gy{\"o}ri, J.~Lehel, and R.~H. Schelp.
\newblock Longest paths in circular arc graphs.
\newblock {\em Combinatorics, Probability and Computing}, 13(3):311--317, 2004.

\bibitem{biro1992precoloring}
M.~Bir{\'o}, M.~Hujter, and Z.~Tuza.
\newblock Precoloring extension. {I}. interval graphs.
\newblock {\em Discrete Mathematics}, 100(1-3):267--279, 1992.

\bibitem{cerioli2020}
M.~R. Cerioli, C.~G. Fernandes, R.~Gómez, J.~Gutiérrez, and P.~T. Lima.
\newblock Transversals of longest paths.
\newblock {\em Discrete Mathematics}, 343(3):111717, 2020.

\bibitem{cerioli2019}
M.~R. Cerioli and P.~T. Lima.
\newblock Intersection of longest paths in graph classes.
\newblock {\em Discrete Applied Mathematics}, 281:96--105, 2020.

\bibitem{chaplick2021kernelization}
S.~Chaplick, F.~V. Fomin, P.~A. Golovach, D.~Knop, and P.~Zeman.
\newblock Kernelization of graph hamiltonicity: Proper {$H$}-graphs.
\newblock {\em SIAM Journal on Discrete Mathematics}, 35(2):840--892, 2021.

\bibitem{chudnovsky2021}
M.~Chudnovsky, J.~King, M.~Pilipczuk, P.~Rz{\k{a}}{\.z}ewski, and S.~Spirkl.
\newblock Finding large {$H$}-colorable subgraphs in hereditary graph classes.
\newblock {\em SIAM Journal on Discrete Mathematics}, 35(4):2357--2386, 2021.

\bibitem{Gallai68}
P.~Erd\H{o}s and G.~Katona, editors.
\newblock {\em Theory of Graphs}.
\newblock Proceedings of the Colloquium held at Tihany, Hungary, September
  1966. Academic Press, New York, 1968.
\newblock Problem 4 (T. Gallai), p. 362.

\bibitem{gao2021nonempty}
Y.~Gao and S.~Shan.
\newblock Nonempty intersection of longest paths in graphs without forbidden
  pairs.
\newblock {\em Discrete Applied Mathematics}, 304:76--83, 2021.

\bibitem{harvey2023intersecting}
D.~J. Harvey and M.~S. Payne.
\newblock Intersecting longest paths in chordal graphs.
\newblock {\em Discrete Mathematics}, 346(4):113284, 2023.

\bibitem{hodur2025finding}
N.~Hodur, M.~Pil{\'s}niak, M.~Prorok, and P.~Rz{\k{a}}{\.z}ewski.
\newblock Finding large {$k$}-colorable induced subgraphs in (bull, chair)-free
  and (bull, {E})-free graphs.
\newblock {\em arXiv preprint arXiv:2504.04984}, 2025.

\bibitem{jobson2016detour}
A.~S. Jobson, A.~E. K{\'e}zdy, J.~Lehel, and S.~C. White.
\newblock Detour trees.
\newblock {\em Discrete Applied Mathematics}, 206:73--80, 2016.

\bibitem{joos2015note}
F.~Joos.
\newblock A note on longest paths in circular arc graphs.
\newblock {\em Discussiones Mathematicae Graph Theory}, 35(3):419--426, 2015.

\bibitem{kang2017width}
D.~Y. Kang, O.-j. Kwon, T.~J. Str{\o}mme, and J.~A. Telle.
\newblock A width parameter useful for chordal and co-comparability graphs.
\newblock {\em Theoretical Computer Science}, 704:1--17, 2017.

\bibitem{kierstead2023improved}
H.~A. Kierstead and E.~R. Ren.
\newblock Improved upper bounds on longest-path and maximal-subdivision
  transversals.
\newblock {\em Discrete Mathematics}, 346(9):113514, 2023.

\bibitem{lima2025longest}
P.~T. Lima and A.~Nikabadi.
\newblock Longest path transversals in claw-free and {$P_5$}-free graphs.
\newblock In {\em International Conference on Algorithms and Complexity}, pages
  310--325. Springer, 2025.

\bibitem{long2021sublinear}
J.~A. Long~Jr, K.~G. Milans, and A.~Munaro.
\newblock Sublinear longest path transversals.
\newblock {\em SIAM Journal on Discrete Mathematics}, 35(3):1673--1677, 2021.

\bibitem{long2023non}
J.~A. Long~Jr, K.~G. Milans, and A.~Munaro.
\newblock Non-empty intersection of longest paths in {$H$}-free graphs.
\newblock {\em arXiv preprint arXiv:2302.07110}, 2023.

\bibitem{long2024longest}
J.~A. Long~Jr, K.~G. Milans, and M.~C. Wigal.
\newblock Longest path and cycle transversals in chordal graphs.
\newblock {\em arXiv preprint arXiv:2412.20729}, 2024.

\bibitem{norin2025small}
S.~Norin, R.~Steiner, S.~Thomass{\'e}, and P.~Wollan.
\newblock Small hitting sets for longest paths and cycles.
\newblock {\em arXiv preprint arXiv:2505.08634}, 2025.

\bibitem{rautenbach2014transversals}
D.~Rautenbach and J.-S. Sereni.
\newblock Transversals of longest paths and cycles.
\newblock {\em SIAM Journal on Discrete Mathematics}, 28(1):335--341, 2014.

\bibitem{walther1969nichtexistenz}
H.~Walther.
\newblock {\"U}ber die nichtexistenz eines knotenpunktes, durch den alle
  l{\"a}ngsten wege eines graphen gehen.
\newblock {\em Journal of Combinatorial Theory}, 6(1):1--6, 1969.

\bibitem{zamfirescu1976longest}
T.~Zamfirescu.
\newblock On longest paths and circuits in graphs.
\newblock {\em Mathematica Scandinavica}, 38(2):211--239, 1976.

\end{thebibliography}

\end{document}